\documentclass[runningheads,french]{asmda2005}
\usepackage[T1]{fontenc}
\usepackage{babel}

\usepackage{asmda2005References} 

\usepackage{graphicx} 

\begin{document}
%
\title*{Missing values : processing with the Kohonen algorithm
}
%
\toctitle{Missing values : processing with the Kohonen algorithm}
%
\titlerunning{Missing values}
%
\author{
  Marie Cottrell
  \and 
  Patrick Letrémy
}
%
\index{Cottrell M.}
\index{Letrémy, P.}

%
\authorrunning{Cottrell and Letrémy}
%
\institute{
  SAMOS-MATISSE\\
  Université Paris 1\\
  90, rue de Tolbiac, 75634 Paris Cedex 13, France\\
  (e-mail: {\tt cottrell@univ-paris1.fr, pley@univ-paris1.fr})
}

\maketitle             

\begin{abstract}
We show how it is possible to use the Kohonen self-organizing algorithm to deal with data with missing values and estimate them. After a methodological reminder, we illustrate our purpose with three applications to real-world data.
Nous montrons comment il est possible d'utiliser l'algorithme d'auto-organisation de Kohonen pour traiter des données avec valeurs manquantes et estimer ces dernières. Après un rappel méthodologique, nous illustrons notre propos à partir de trois applications à des données réelles.
\end{abstract}

\section{Introduction}


The processing of data which contain missing values is a complicated and always awkward problem, when the data come from real-world contexts. In applications, we are very often in front of observations for which all the values are not available, and this can occur for many reasons: typing errors, fields left unanswered in surveys, etc.

Most of the statistical software (as SAS for example) simply suppresses incomplete observations. It has no practical consequence when the data are very numerous. But if the number of remaining data is too small, it can remove all significance to the results.

To avoid suppressing data in that way, it is possible to replace a missing value with the mean value of the corresponding variable, but this approximation can be very bad when the variable has a large variance.

So it is very worthwhile seeing that the Kohonen algorithm (as well as the Forgy algorithm) perfectly deals with data with missing values, without having to estimate them beforehand. We are particularly interested in the Kohonen algorithm for its visualization properties.

In Smaïl Ibbou's PHD thesis, one can find a chapter about this question, but it has not been published yet. The examples are run with the software written by Patrick Letrémy in IML-SAS and available on the SAMOS WEB page (http://samos.univ-paris1.fr).

\section{Adaptation of the Kohonen algorithm to data with missing values}

We do not remind of the definition of the Kohonen algorithm here, see for example Kohonen \cite{kohonen}, or \cite{cottrell}.

Let us assume that the observations are real-valued $p$-dimensional vectors, that we intend to cluster into $n$ classes.

When the input is an incomplete vector $x$, we first define the set $M_x$ of the numbers of the missing components. $M_x$ is a sub-set of $\{1, 2, \ldots, p\}$. If $C=(C_1, C_2, ..., C_n)$ is the set of code-vectors at this stage, the winning code-vector $C_{i_0(x,C)}$ related to $x$ is computed as by setting 
$$ i_0(x,C)=\arg \min_i \|x-C_i\|,$$

where the distance $\|x - C_i\|^2 = \sum_{k \not\in M_x}(x_k - C_{i,k})^2$ 
is computed with the components present in vector $x$.

One can use incomplete data in two ways: 

a) If we want to use them during the construction of the code-vectors, at each stage, the update of the code-vectors (the winning one and its neighbors) only concerns the components present in the observation. Let us denote 
$C^t=(C_1^t, C_2^t, ..., C_n^t)$ the code-vectors at time $t$ and if a randomly chosen observation $x^{t+1}$ is drawn, the code-vectors are updated by setting: 

$$C_{i,k}^{t+1}= C_{i,k}^t + \epsilon (t) (x_k^{t+1} - C_{i,k}^t)$$
for $k\notin M_x$ and $j$ neighbor of $i_0(x^{t+1}, C^t)$. Otherwise, 
$$C_{i,k}^{t+1}= C_{i,k}^t.$$

The sequence $\epsilon(t)$ is [0,1]-valued with $\epsilon(0)\simeq 0.5$ and converges to 0 as $1/t$. After convergence, the classes are defined by the nearest neighbor method.

b) If the data are numerous enough to avoid using the incomplete vectors to build the map, one can content oneself with classifying them after the map is built, as supplementary data, by allocating them to the class with the code-vector which is the nearest for the distance restricted to non-missing components.

This method yields excellent results, provided a variable is not totally or almost totally missing, and also provided the variables are correlated enough, which is the case for most real data bases. Several examples can be encountered in Smaïl Ibbou's PHD thesis \cite{ibbou} and also in Gaubert, Ibbou and Tutin \cite{gaubert}.

\section{Estimation of missing values, computation of membership probabilities}

Whatever the method used to deal with missing values, one of the most interesting properties of the algorithm is that it allows an a posteriori estimation of these missing values.

Let us denote by $C=(C_1, C_2, \ldots, C_n)$ the code-vectors after building the Kohonen map. If $M_x$ is the set of missing component numbers for the observation $x$, and if $x$ is classified in class $i$, for each index $k$ in $M_x$, one estimates $x_k$ by : 
$$\hat{x}_k=C_{i,k}.$$

Because in the end of the learning the Kohonen algorithm uses no more neighbor (0 neighbor algorithm), we know that the code-vectors are asymptotically near the mean values of their classes. This estimation method therefore consists in estimating the missing values of a variable by the mean value of its class.

It is clear that this estimation is all the more precise as the classes built by the algorithm are homogeneous and well separated. Numerous simulations have shown as well for artificial data as for real ones, that when the variables are sufficiently correlated, the precision of these estimations is remarkable, \cite{ibbou}.

It is also possible to use a probabilistic classification rule, by computing the membership probabilities for the supplementary observations (be they complete or incomplete), by putting:
$$Prob(x \in \mbox{Class }i) = 
\frac{\exp (-\|x - C_i\|^2)}{\sum^n_{k=1}\exp (-\|x - C_k\|^2)}.$$

These probabilities also give confirmation of the quality of the organization in the Kohonen map, since significant probabilities have to correspond to neighboring classes.

Moreover, to estimate the missing values, one can compute the weighted mean value of the corresponding components. The weights are the membership probabilities. If $x$ is an incomplete observation, and for each index $k$ in $M_x$, one estimates $x_k$ by : 
$$\hat{x}_k=\sum Prob(x \in \mbox{Class }i) \; C_{i,k}.$$

These probabilities also provide confidence intervals, etc. In the following sections, we present three examples extracted from real data.

\section{Socio-economic data}

The first example is classical. The database contains seven ratios measured in 1996 on the macroeconomic situation of 182 countries. This data set was first used by F. Blayo and P. Demartines \cite{blayo} in the context of data analysis by SOMs. 

The measured variables are: annual population growth (ANCRX), mortality rate (TXMORT), illiteracy rate (TXANAL), population proportion in high school (SCOL2), GDP per head (PNBH), unemployment rate (CHOMAG), inflation rate (INFLAT).

Among the set of 182 countries, only 115 have no missing values, 51 have only one missing value, while 16 have 2 or more than 2 missing values.

Therefore we use the 115 + 51 = 166 complete or almost complete countries to build the Kohonen map, and we then classify the 16 remaining countries. The data are centered and reduced as classically. We take a Kohonen map with 7 by 7 units, that is 49 classes. Figure 1 shows the contents of the classes. The 166 countries that were used for computing the code-vectors are in normal font, the 16 others in underlined italics. 

\begin{figure}
\caption{The 182 countries (166 + 16) on a 7 by 7 map, 1500 iterations}
\includegraphics[scale=0.5]{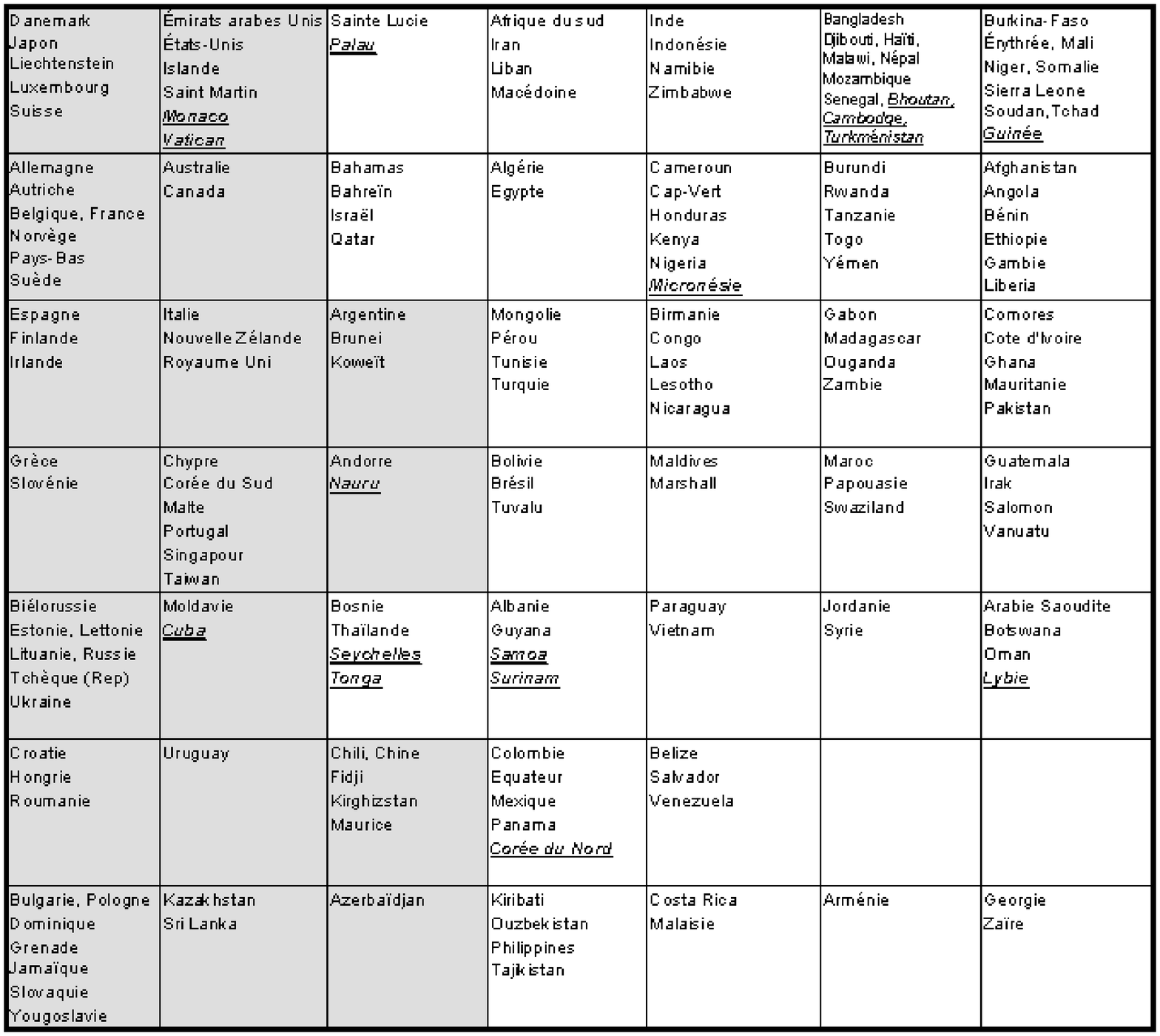}
\end{figure}

We can see that rich countries are in the top left hand corner, very poor ones are displayed in the top right hand corner. Ex-socialist countries are not very far from the richest, etc. As for the 16 countries which are classified after the learning as supplementary observations, we observe that the logic is respected. Monaco and Vatican are displayed with rich countries, and Guinea with very poor countries, etc.

From these computations, it is possible to calculate the membership probabilities of each supplementary observation of each of the 49 classes.

For example, the probabilities that Cuba belongs to class $i$ are greater than 0.03 for classes $i= (1,1), (2,1), (3,1), (4,1), (5,1), (6,1), (7, 1), (1,2), (2,2)$,
$ (3,2), (4,2), (5,2), (6,2), (7,2), (3,3), (4,3), (6,3), (7,3)$, the maximum (0.06) being reached for class $(5,2)$. We can notice (figure 1) that they are neighboring classes. From these probabilities, it is possible to estimate the distribution of the estimators of the missing values. For Cuba, the variables in question are GDP, Unemployment and Inflation.

From these results, it is possible (as it will be shown in the talk) to build super-classes by using an ascending hierarchical classification of the code-vectors and then to cross this classification with other exogenous classifications, etc.

\section{Study of the property market in Ile-de-France}

The second example is extracted from a study commissioned by the direction of Housing in the Regional Direction of Equipment in Ile-de-France (DHV/DREIF). This was achieved in 1993 by Paris 1 METIS and SAMOS laboratories, by Gaubert, Tutin and Ibbou, \cite{gaubert}.

For 205 towns in Ile-de-France considered in 1988, we have property data (housing rents and prices, old and new, collective or individual, standard or luxurious, office rents and prices, old and new). Structurally, some of the data are missing, for example the office market can be nonexistent in some towns.

This is a case where some data are structurally missing, and where the number of towns is dramatically reduced if one suppresses those which are incomplete: only 5 out of 205 would be kept! So for the learning, we use 150 towns which have less than 12 missing values out of 15. After that, the 55 towns which have more than 12 missing values out of 15 are classified as supplementary observations.

Figure 2 displays the 205 towns (with and without missing values) classified on a 7 by 7 Kohonen map. Note that there are about 63\% of missing values on the data set.

\begin{figure}
\caption{The 205 towns in Ile-de-France, in underlined italics the 55 towns which have more than 12 missing values out of 15 variables}
\includegraphics[scale=0.5]{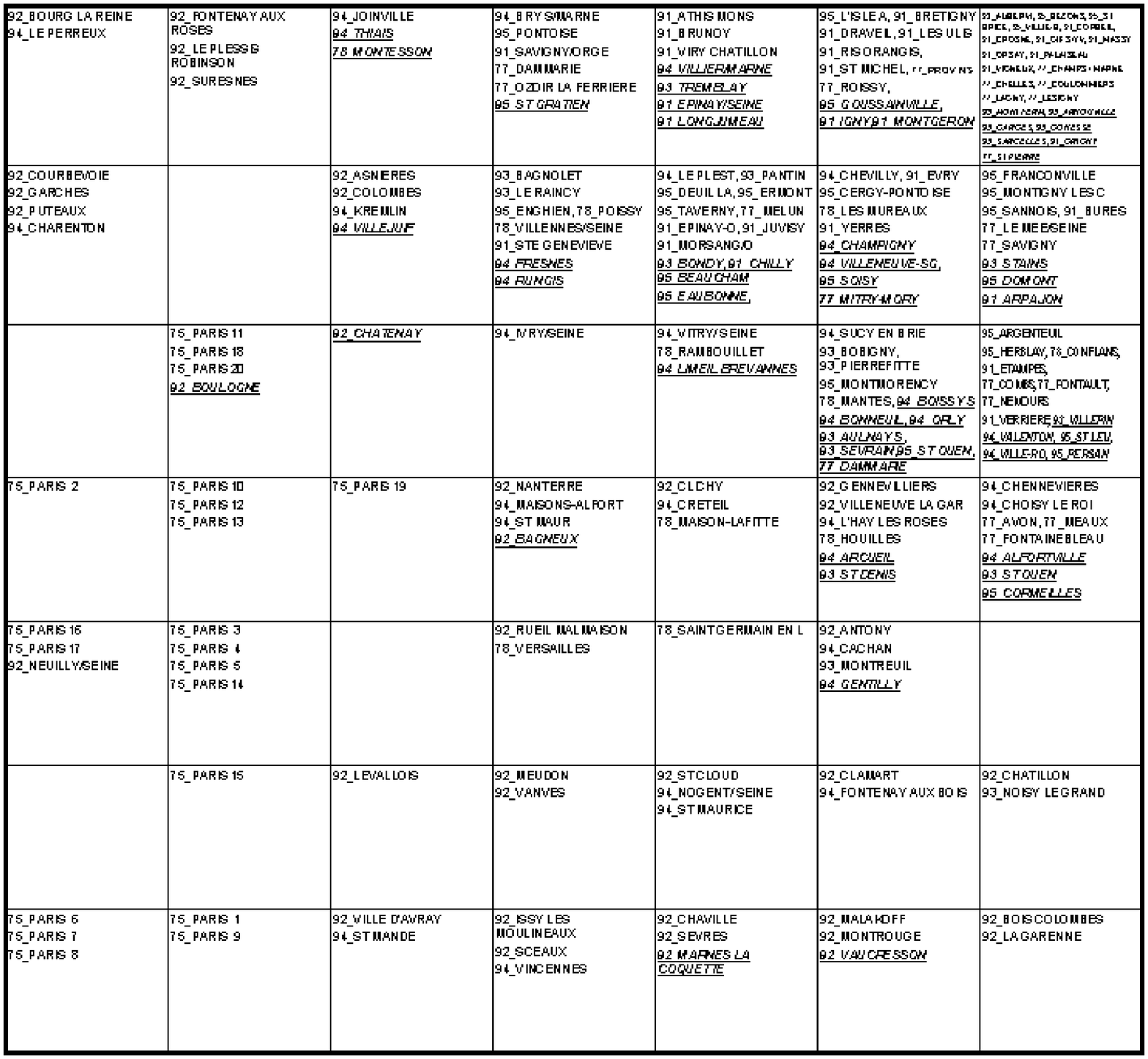}
\end{figure}

In this example which is practically impossible to deal with using classical software, we see that the Kohonen algorithm nevertheless allows to classify extremely sparse data, without introducing any rough error. 
The results are perfectly coherent, even though the data are seriously incomplete. The districts of Paris, Boulogne and Neuilly sur Seine are in the bottom left hand corner. On a diagonal stripe, one finds the towns of the inner suburbs (petite couronne), further right there are the towns of the outer suburbs (grande couronne). Arcueil is classified together with l'Ha\"y-les-Roses (class (2,3)), Villejuif with Kremlin Bicêtre (4,6), etc.

Of course, these good results can be explained by the fact that the 15 measured variables are well correlated and that the present values contain information about missing values. The examination of the correlation matrix (that SAS computes even in case of missing values) shows that 76 coefficients out of 105 are greater than 0.8, none of them being less than 0.65.

\section{Structures of Government Spending from 1872 to 1971}

The third example is a very classical one in data analysis, taken from the book ``Que-sais-je ?'' by Bouroche and Saporta , ``L'analyse des données''  \cite{saporta}. The problem is to study the government spending, measured over 24 years between 1872 and 1971, by a 11-dimensional vector: Public Authorities (Pouvoirs publics), Agriculture (Agriculture), Trade and Industry (Commerce et industrie), Transports (Transports), Housing and Regional Development (Logement et aménagement du territoire), Education and Culture (Education et culture), Social Welfare (Action sociale), Veterans (Anciens combattants), Defense (Défense), Debt (Dette), Miscellaneous (Divers). It is a very small example, with 24 observations of dimension 11, without any missing values.

A Principal Component Analysis provides an excellent representation in two dimensions with 64\% of explained variance. See figure 3.

\begin{figure}
\caption{On the left, the projections on the first two principal axes; on the right, the Kohonen map with 9 classes and 3 super-classes}
\includegraphics[scale=0.5]{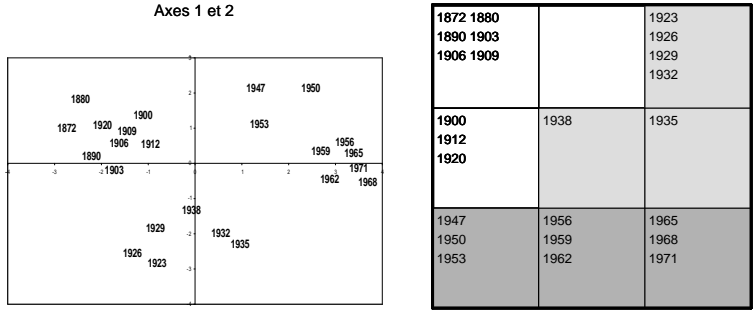}
\end{figure}

On this projection, the years split up into three groups, which correspond to three clearly identified periods (before the First World War, between the two World Wars, after the Second World War). Only the year 1920, the first year when an expenditure item for Veterans appears, is set inside the first group, while it belongs to the second one. On the Kohonen map, the three super-classes (identical to the ones just defined) are identified by an ascending hierarchical classification of the code-vectors.

In this example, we have artificially suppressed randomly chosen values which were present in the original data, from 1 value out of 11 to 8 values out of 11, in order to study the clustering stability and compute the accuracy of the estimations that we get by taking the corresponding values of the code-vectors.

One can observe that the three super-classes remain perfectly stable as long as one does not suppress more than 3 values a year, that is 27\% of the values.

Then we estimate the suppressed values in each case. The next table shows the evolution of the mean quadratic error according to the number of suppressed values.

\begin{table}
  \begin{center}
    \begin{tabular}{|l|l|l|l|l|l|l|l|l|}
      \hline
      Number of missing values & 1 & 2 & 3 & 4 & 5 & 6 & 7 & 8\\
		 \hline
Percentage of missing values & 9\% & 18\% & 27\% & 36\% & 45\% & 55\% & 64\% & 73\%\\
      \hline
      & 0.39 & 0.54 & 0.73 & 1.11 & 1.31 & 1.30 & 1.27 & 1.39\\
           \hline
    \end{tabular}
  \end{center}
  \caption{Mean Quadratic Error according to the number of suppressed values}
 
\end{table}

We notice that the error remains small as long as we do not suppress more than 3 values a year. 

\section{Conclusion}

Through these three examples, we have thus shown how it is possible and desirable to use Kohonen maps when the available observations have missing values. Of course, the estimations and the classes that we get are all the more relevant since the variables are well correlated.

Example 2 shows that it can be the only possible method when the data are extremely sparse. Example 3 shows how this method allows to estimate the absent values with good accuracy. The completed data can then be dealt with using any classical treatment.

\bibliographystyle{asmda2005References}

\bibliography{bibliography-cottrell}

\end{document}